\documentclass[12pt]{article}
\usepackage{latexsym, amssymb}
\usepackage[all]{xy}
\def\separation{\medskip}

\def\Q{{\bf Q}}
\def\P{I\!\!P}

\def\k{{\bf k}}

\def\L{{\cal L}}

\def\I{{\cal I}}
\def\D{{\cal D}}

\def\?{{\bf ??}}

\def\rig#1{\smash{ \mathop{\longrightarrow}\limits^{#1}}}

 \newtheorem{theorem}{Theorem}[section]

\newtheorem{prop}[theorem]{Proposition}
\newtheorem{definition}[theorem]{Definition}
\newtheorem{corollary}[theorem]{Corollary}

\newtheorem{remark}[theorem]{Remark}

\newcommand{\proof}{{\it Proof.}\ }
\newcommand{\qed}{\hfill  $\Box$\separation}

\begin{document}

\title{On singular L\"uroth quartics}
 \author{Giorgio Ottaviani - Edoardo Sernesi\footnote{Both authors are members of GNSAGA-INDAM.}}
\date{}
 \maketitle

\begin{small}\emph{  Dedicated to Fabrizio Catanese on the occasion of his
60-th birthday}\end{small}

 \abstract{Plane quartics containing  the ten
  vertices of a complete pentalateral and limits of them are called L\"uroth quartics.
The locus of singular L\"uroth quartics has two  irreducible
components, both  of codimension two in $\P^{14}$. We compute the
degree of them and we discuss the consequences of this computation
on the explicit form of the L\"uroth invariant. One important tool
are the Cremona hexahedral equations of the cubic surface. We also
compute the class in $\overline{M}_3$ of the closure of the locus
of nonsingular L\"uroth quartics.}

\section*{Introduction}

All schemes and varieties will be assumed to be defined   over an
algebraically closed field $\k$ of characteristic zero.
 We recall that a \emph{complete pentalateral} in $\P^2$
is a configuration consisting of five lines, three by three
linearly independent, together with the ten double points of their
union, which are called \emph{vertices} of the pentalateral.
 A \emph{nonsingular
L\"uroth quartic} is a nonsingular quartic plane curve containing
the ten vertices of a complete pentalateral. Such curves fill an
open set of an irreducible, SL$(3)$-invariant, hypersurface
$\L\subset\P^{14}$. The (possibly singular) quartic curves
parametrized by the points of $\L$ will be called \emph{L\"uroth
quartics}. In   \cite{OS09} we have computed that $\L$ has degree
$54$, by reconstructing a proof published by Morley in 1919
\cite{fM14}. Another proof has been given by Le Potier and
Tikhomirov in \cite{PT}.

In this paper we put together the projective  techniques of
\cite{fM14} and \cite{OS09} with the cohomological techniques in
\cite{PT}, and we prove some new results about the L\"uroth
hypersurface.
 We refer to the introduction of \cite{OS09} for an explanation of
  the connection of this topic with moduli of vector bundles on $\P^2$.

The locus of singular L\"uroth quartics has been considered in
\cite{fM14} and \cite{PT}. It is obtained as the intersection
between the L\"uroth hypersurface of degree $54$ and the
discriminant of degree $27$. It has two irreducible components
$\L_1$ and $\L_2$, both  of codimension $2$ in $\P^{14}$,  and it
is known that $\deg (\L_2)_{red}=27\cdot 15$ (\cite{PT}, Cor.
9.4). We compute the degree of $\L_1$, a question left open in
\cite{PT} (end of 9.2). Indeed we prove the following theorem.

\begin{theorem}\label{T:main0}
 The intersection between the L\"uroth hypersurface $\L$
 and the discriminant $\D$ is transverse along $\L_1$, and

(i) $\mathrm{deg} \L_1=27\cdot 24$.

(ii) $\L_2$ is non reduced of degree $27\cdot 30$,
 \end{theorem}

We warn the reader that our $(\L_2)_{red}$ corresponds to $\L_2$ in \cite{PT}.

An interesting aspect of the geometrical construction in
\cite{OS09} is that any smooth cubic surface $S$ defines in a
natural way $36$ planes, which we called Cremona  planes, one for
each  of the $36$ double-six configurations of lines on $S$. Their
main property is that the ramification locus of the projection
$\pi_P$ centered at $p\in S$ is a L\"uroth quartic if and only if
$p$ belongs to any of the Cremona  planes.

We describe the Cremona  planes on a  nonsingular cubic surface by
pure projective geometry. To give the flavour of this construction
we state the following result.

\begin{theorem}\label{noninvolutory}
Let $S$ be a nonsingular cubic surface. Fix a double-six on $S$.
Let $\ell_s$,  $s=1,\ldots ,15$,  be the $15$ remaining lines. For
each $1\le s \le 15$ consider the three planes $\Pi_{s,h}$,
$h=1,2,3$, containing $\ell_s$ such that $S\cap\Pi_{s,h}$ consists
of $\ell_s$ and of two residual lines not belonging to the
double-six. Let $P_{s,h}$ be the intersection point of the two
residual lines. Then the  $15$ points $\ell_s\ \cap <P_{s,1},
P_{s,2}, P_{s,3}>$, $s=1, \dots, 15$, lie on a plane, which is the
Cremona plane associated to the double-six.

\end{theorem}

  Theorem \ref{noninvolutory2}
   contains the statement of this theorem with
additional informations on the {\it involutory} and {\it non
involutory} points. They give a geometrical explanation of the
reducibility of $\L\cap \D$, see Prop. \ref{geominvolutory}.

 We conclude  with the statement of non-existence
 of an invariant of degree 15 (Prop. \ref{noinv})
vanishing on $(\L_2)_{red}$,  which we have obtained by a computer
computation. This means that $(\L_2)_{red}$ is not a complete
intersection,  and we relate this fact with the last sentence in
Morley's paper \cite{fM14}. This leads to a reconstruction of some
speculations of Morley about the (still unknown) explicit form of
the L\"uroth hypersurface. Our result implies that these
speculations are partially wrong, but with a slight correction
they might become   true.

In 1967 Shioda \cite{Sh} found the Hilbert series for the
invariant  ring of plane quartics. From his formula it follows
that the space of invariants of degree $54$ has   dimension
$1165$. This shows the difficulty to find the explicit expression
(or the symbolic expression) of the L\"uroth invariant, which, to
the best of our knowledge, is still unknown. In the last section
we compute the class in $\overline{M}_3$ of the divisor of
L\"uroth quartics.

The content of the paper is the following.

\noindent In the first section we summarize some results from \cite{OS09}
on Cremona hexaedral equations and Cremona  planes on a cubic surface. We recall the
 purely geometric construction of the involutory points.

\noindent In the second  section we summarize well known facts
about the description of plane quartics as symmetric determinants
with linear entries. We recall how L\"uroth quartics can be found
in this description (they are the image of a pfaffian hypersurface
$\Lambda$,
 which is an invariant
 of degree $6$) and, following
\cite{PT} and \cite{Giz}, we also  describe how the two components of singular
L\"uroth quartics can be found.

\noindent The third section contains our new results, the main
ones being described above, and their proofs. The last section is
devoted to the computation  the class $[L]$ of the divisor in
$\overline{M}_3$ parametrizing L\"uroth quartics, as well as to
some related remarks.

We thank I. Dolgachev for calling to our attention the reference
\cite{Giz} and C. Faber for a helpful  conversation with the
second author about the topics of the last section.

    \section{Cremona hexahedral equations and Cremona  planes}\label{S:cremona}

\noindent
    Recall that a  \emph{double-six} of lines on
    a nonsingular cubic surface $S \subset \P^3$ consists of two sets of six lines
    $\Delta = (A_1, \dots, A_6; B_1, \dots B_6)$
    such that the lines $A_j$ are mutually skew as
     well as the lines  $B_j$; moreover each $A_i$ meets each $B_j$ except when $i=j$.

   \noindent
   In  $\P^5$ with coordinates $(Z_0, \dots,Z_5)$  consider the following equations:
    \begin{equation}\label{E:cre1}
    \cases{Z_0^3+Z_1^3+Z_2^3+Z_3^3+Z_4^3+Z_5^3 = 0 \cr
    Z_0+Z_1+Z_2+Z_3+Z_4+Z_5 = 0  \cr
    \beta_0 Z_0+\beta_1 Z_1 + \beta_2 Z_2 + \beta_3 Z_3 + \beta_4 Z_4 + \beta_5 Z_5 = 0 }
    \end{equation}
    where the $\beta_s$'s are general
     constants. These equations define a nonsingular cubic  surface $S$ in a $\P^3$ contained
    in $\P^5$ and are called \emph{Cremona hexahedral equations} of $S$, after \cite{Cre}.

For any choice of two disjoint pairs  of indices
$\{i,j\}\cup\{k,l\}\subset\{0,\ldots 5\}$, the equations
$Z_i+Z_j=Z_k+Z_l=0$ define a line contained in $S$. There are $15$
such lines and the remaining $12$  determine a double-six of lines
on $S$. Therefore the equations (\ref{E:cre1}) define a double six
on $S$.  More precisely we have the following:

    \begin{theorem}\label{T:cre1}
    Each system of Cremona hexahedral equations
     of a nonsingular cubic surface $S$ defines a double-six
      of lines on $S$. Conversely, the choice of a double-six
       of lines on $S$ defines a   system of Cremona hexahedral equations
       (\ref{E:cre1}) of $S$, which is uniquely determined up to replacing
       the coefficients $(\beta_0, \dots, \beta_5)$ by $(a+b\beta_0,\dots, a+b\beta_5)$ for some $a,b \in \k$, $b \ne 0$.
    \end{theorem}

    We refer to \cite {D}, Theorem  9.4.6,  for the proof.   We need to point out from \cite{OS09}, Coroll. 4.2 the following:

     \begin{corollary}\label{C:cre1}
     To a pair $(S,\Delta)$ consisting of a
     nonsingular cubic surface $S \subset \P^3$ and a  double-six of lines $\Delta$ on $S$,
     there is canonically associated a plane $\Xi \subset \P^3$   which is given by the equations
          \begin{equation}\label{cre4}
     \cases{
    Z_0+Z_1+Z_2+Z_3+Z_4+Z_5 = 0  \cr
    \beta_0 Z_0+\beta_1 Z_1 + \beta_2 Z_2 + \beta_3 Z_3 + \beta_4 Z_4 + \beta_5 Z_5 = 0 \cr
    \beta_0^2Z_0+\beta_1^2Z_1 + \beta_2^2 Z_2 +\beta_3^2Z_3 +\beta_4^2Z_4+  \beta_5^2Z_5 = 0}
     \end{equation}
     where the coefficients $\beta_0,\dots,\beta_5$ are those
     appearing in the Cremona equations of $(S,\Delta)$.
    \end{corollary}

      \begin{definition}\label{D:cocre}
    The plane $\Xi \subset \P^3$ will be called the \emph{Cremona  plane} associated to the pair  $(S,\Delta)$.
    \end{definition}

 The link with the L\"uroth quartics is given by the following.

\begin{theorem}\label{T:cocrelur} Let $p\in S$ be a point.
 The projection from $p$ defines a rational double
covering $\xymatrix{\pi_P\colon S\ar@{-->}[r]& \P^2}$  ramified over a plane quartic.
The ramification curve is a L\"uroth quartic if and only if $p$ belongs to any
of the Cremona  planes.
\end{theorem}
\proof
See the Remark 10.7 and Theorem 6.1 of \cite{OS09}.
 \qed

There is a second description of the Cremona  planes, by means of the involutory points.
In order to state it, we give, following \cite{OS09}, a geometric construction  of the
involutory points.
Consider two skew lines $A,B \subset S$.   Denote by
   $ f\colon A\rightarrow B$
   the double cover associating to $p \in A$ the point
    $f(p):= T_pS \cap B$ where $T_pS$ is the tangent plane to
   $S$ at $p$.
  Define
  $ g\colon B \rightarrow A$
 similarly.
Let $p_1,p_2 \in A$  (resp. $q_1,q_2 \in B$) be the ramification points of $f$ (resp. $g$).
Consider the pairs of branch points
$f(p_1),f(p_2) \in B,  \  g(q_1), g(q_2) \in A$,
and the new  morphisms
\[
f'\colon A \rightarrow \P^1, \qquad g'\colon B \rightarrow \P^1
\]
defined by the conditions that  $g(q_1), g(q_2)$ are  ramification
points of $f'$ and $f(p_1),f(p_2)$ are ramification points of
$g'$. Let   $Q_1+Q_2$ (resp. $P_1+P_2$) be the   common divisor of
the two $g^1_2$'s on $A$ (resp. on $B$) defined by $f$ and $f'$
(resp. by  $g$ and $g'$).
  The points
 \[
 \bar P= g(P_1)=g(P_2) \in A, \qquad  \bar Q =f(Q_1)=f(Q_2)\in B
 \]
 are called  the \emph{involutory points} (relative to the pair of lines $A$ and $B$).

Note that each line $A \subset S$ contains 16 involutory points,
which correspond to the 16 lines $B \subset S$ which are skew with
$A$, and they are distinct (see the proof of Prop. 6.3 of
\cite{OS09}).

Let $\bar P_i \in A_i, \  \bar Q_i \in B_i$ be the   involutory
points  relative to the pair $A_i$ and $B_i$. We obtain twelve
points
\[
\bar P_1, \dots, \bar P_6, \bar Q_1, \dots, \bar Q_6 \in S
\]
which are canonically associated to the double-six  $\Delta$.

\begin{theorem}\label{T:plane1}
For any double-six $\Delta$ there is a unique plane  \  $\Xi \subset \P^3$ containing  the involutory points
\[
\bar P_1, \dots, \bar P_6, \bar Q_1, \dots, \bar Q_6
\]
   Moreover $\Xi$ coincides with the Cremona  plane associated to the pair $(S,\Delta)$.

The $36$ Cremona  planes obtained in this way are distinct.
\end{theorem}

\proof See \cite{OS09} Theorem 6.1 and Prop. 6.3.
\qed

\section{The symmetric representation of L\"uroth
quartics}\label{S:symm}

Let $Q_0$, $Q_1$, $Q_2$ be three  linearly  independent quadrics
in $\P^3={\P}(W)$. They generate a net of quadrics $\langle
Q_0,Q_1,Q_2\rangle$ whose base locus, in general, consists of
eight points in general position. We can parametrize any net of
quadrics by the points of $\P^2=\P(V)$, and as such it can be seen
as an element $f\in{\P}(V\otimes S^2W)$.  The symmetric
determinantal representation of the quadrics of the net gives a
dominant rational  map $ \delta\colon {\P}(V\otimes
S^2W)\dashrightarrow {\P}(S^4V)$, see \cite{D}.
 In \cite{W} Wall studied  the map $\delta$ in the setting of invariant theory. He proved
that the non-semistable points for the action of $SL(W)$ on ${\P}(V\otimes S^2W)$
are exactly given by the locus $Z(\delta)$ where $\delta$ is not defined.

There is a factorization
 through the GIT quotient

$$\xymatrix{
\P{{}}(V\otimes S^2W)^{ss}\ar[d]^\pi \ar[dr]^\delta \\
\P{{}}(V\otimes S^2W)//SL(W)\ar[r]^-g & \P{{}}(S^4V)}
 $$
where $g$ is generically finite of degree $36$  and $\P(V\otimes
S^2W)//SL(W)$ parametrizes pairs $(B,t)$ consisting of a plane
quartic $B$ and an even theta-characteristic $t$ on it.

  \noindent
 Consider the hypersurface $\Lambda\subset \P(V\otimes S^2W)$ of
 degree 6 consisting of the nets $\langle Q_0,Q_1,Q_2\rangle$ satisfying the equation:
 $$\textrm{Pf}\left[\begin{array}{rrr}
0&Q_0&-Q_1\\
-Q_0&0&Q_2\\
Q_1&-Q_2&0\\
\end{array}\right]=0$$
where we identify each quadric $Q_i$ with its  corresponding
symmetric matrix.  It can be shown \cite{hB14} that a net belongs
to $\Lambda$ if and only if  $Q_0,Q_1,Q_2$ are the polar quadrics
of three points with respect to a cubic surface in $\P^3$.

\noindent
In \cite{gO07}, section 4, it is shown that $\Lambda$ is
the   $5$-secant variety  of ${\P}(V)\times{\P}(W)$ embedded with $ {\cal O}(1,2)$,
and this fact is used to give a new  proof of the L\"uroth theorem.

Let $\L \subset S^4(V)$ be the L\"uroth invariant  of degree $54$.
Since the entries of $\delta$ have degree four,   $\delta^*\L$ has
degree $216$ and the crucial fact, for our purposes, is that
$\delta^*\L$  contains $\Lambda$ as an irreducible component (see
for example \cite{gO07}, prop. 6.3(ii)).

\noindent Then $g^{-1}(\L)$ decomposes into two irreducible
components $P$ and $\widetilde P$, both dominating $\L$, with
degree $1$ and $35$ respectively, and $P$ generically parametrizes
the pairs $(B,t)$ where  $t$ is the pentalateral
theta-characteristic on $B$ (see Remark 10.7 of \cite{OS09}). In
particular, $g$ has a rational section over $\L$.

There are two other  classically known invariants  of nets of
quadrics with respect to $SL(V)\times SL(W)$, nicely reviewed by
Gizatullin in \cite{Giz}.
 The \emph{tact-invariant} $J$ has degree $48$ and  vanishes if and only if two of the eight base points
of the net coincide.
The invariant $I$ of degree $30$ vanishes when the net contains a quadric of rank $\le 2$.

\noindent Let $\D$ be the discriminant invariant, which  is
irreducible of degree $27$;  then $\delta^*\D$ is an invariant  of
degree $108$ of the nets of quadrics. Salmon proved the beautiful
identity (up to scalar constants)
\begin{equation}\label{salmonid}
\delta^*\D=I^2J
\end{equation}
It can be interpreted as saying  that  there are two ways to get a
singular quartic as a symmetric determinant. This is interesting
when applied to L\"uroth quartics. The singular L\"uroth quartics
are the elements of  $\L\cap\D$. In \cite{PT}, \S 9, it is shown
that  this locus has two irreducible components $\L_1$ and $\L_2$,
so that we have necessarily
$$\L_1=\delta(\Lambda \cap \{J=0\})\qquad (\L_2)_{red}=\delta(\Lambda\cap \{I=0\})$$
(see Proposition \ref{geominvolutory}). To connect our description
with the setting of \cite{PT} it is enough to note that the
geometric quotient of $\Lambda$ by $SL(W)$ is isomorphic to a
compactification $P$ of the moduli space $M(0,4)$ of rank $2$
stable bundles on $\P^{2}$ with $c_1=0$ and $c_2=4$ and the
restriction of $g$ to $M$ can be identified with the Barth map
(\cite{gO07}, \S 8).

\noindent The fact that the locus of singular L\"uroth quartics
consists of two irreducible components was known also to Morley.

\noindent The $36$ elements of the $g$-fiber over  a general point
of $\D$ are of two types:  there are  $16$ points in
$\pi(\{J=0\})$ and $10$ double points in $\pi(\{I=0\})$.
 This decomposition corresponds to the
  two types of even theta-characteristics on a quartic nodal
  curve:
  16 of them are represented by invertible sheaves, and 10 by torsion-free
   non-invertible sheaves, each counted with multiplicity two (\cite{Giz}, remark 10.1, and \cite{Har}).

\section{The main results and their proofs}\label{S:concluding}

Consider the projective bundle $\pi\colon \P(\Q) \rightarrow
\P^3$ where $\Q= T_{\P^3}(-1)$ is the tautological quotient
bundle. For each $z \in \P^3$ the fibre $\pi^{-1}(z)$ is the
projective plane of lines through $z$. Also consider the
projective bundle $ \beta\colon \P(S^4\Q^\vee) \rightarrow \P^3$.
  For each $z \in \P^3$ the fibre  $\beta^{-1}(z)$ is
  the linear system of quartics in  $\pi^{-1}(z)$.
The Picard group of $\P=\P(S^4\Q^\vee)$ is generated by $H= {\cal O}_{\P}(1)$ and by
  the pullback  $F$ of a  plane in $\P^3$.
  Let
$ \widetilde\L \subset \P(S^4\Q^\vee)$ be the $\beta$-relative
hypersurface  of L\"uroth quartics.  It is invariant under the
natural action of SL$(4)$ on $\P(S^4\Q^\vee)$,  and in \cite{OS09}
we showed that $ \widetilde\L =54H-72F$. Moreover every invariant
of a plane quartic of degree $d$ gives a covariant of the cubic
surface of degree $\frac{2d}{3}$, see \cite{OS09}, Remark 8.2.

\noindent Consider  also the relative invariant subvarieties
$\widetilde \L_1$, $\widetilde\L_2$ and $\widetilde \D$ in the
projective bundle $\P(S^4{\bf Q}^{\vee})$ on $\P^3$. For every
smooth cubic surface $S$ we have the projection
$\P(S^4\Q^\vee)_{|S}\rig{\beta} S$ and have defined the section
$s: S \to \P(S^4\Q^\vee)_{|S}$  associating to $p\in S$ the branch
curve of the projection from $p$.
 It is well known that $s^*(\widetilde \D)$
 consists of the divisor of the twenty-seven lines,
with multiplicity two, cut indeed by a covariant of $S$ of degree $18$.

\noindent In   Theorem \ref{T:cocrelur} we have proved that
$s^*(\widetilde \D\cap \widetilde \L)$ consists of the
intersection of the divisor $s^*(\widetilde \D)$ with the $36$
 Cremona planes.
This is a zero dimensional scheme, consisting of two parts: the involutory points (see Theorem \ref{T:plane1}) and
the non-involutory points. By the above, its length is given by
$$
\deg \D\cdot\deg \L\cdot\left(\frac{2}{3}\right)^2\cdot 3=27\cdot 54\cdot\frac{4}{3}=27\cdot 72
$$
 Therefore on every line on $S$ the
 scheme $s^*(\widetilde \D\cap \widetilde \L)$ has length $72$,
 and multiplicity $\ge 2$ at each point.  It is supported on the $16$ involutory points and on
$\le 20$ non-involutory points.

 \noindent

\begin{prop}\label{geominvolutory}
(i) projecting $S$ from an involutory point of $s^*(\widetilde
\D\cap \widetilde \L)$ we get a branch quartic in $\L_1$
(corresponding to the tact-invariant $J$).

(ii) projecting $S$ from a non-involutory point of $s^*(\widetilde
\D\cap \widetilde \L)$  we get a branch quartic in $\L_2$
(corresponding to the invariant $I$).

(iii) On each line of $S$ there are exactly 10 non-involutory
points, each counts with multiplicity four in $s^*(\widetilde
\D\cap \widetilde \L)$.
\end{prop}
\proof By the Salmon identity (\ref{salmonid}), the two types of
points correspond to the vanishing of the two invariants $I$ and
$J$. We have just to distinguish which is the  type obtained by
each invariant. A check on the degrees
($\frac{60}{48}=\frac{20}{16}$) suffices to prove (i) and (ii). In
order to prove (iii), consider that in the Salmon identity
(\ref{salmonid}) the invariant $I$ appears with exponent two, and
this implies that the non-involutory
 points have to be double ones, that is there are $10$ distinct non-involutory points on each
line,
and at each of these points two Cremona  planes meet.
\qed

  Proposition \ref{geominvolutory} explains
  why $\L_2$ is non reduced. We will construct directly the ten non-involutory points
in   Theorem \ref{noninvolutory2}.

\begin{remark}\label{R:hunt}\rm
It appears that the $36$ Cremona  planes  carry an
interesting combinatorial configuration. Each of them has 27
marked points given by the intersection with the lines on the
cubic surface, 12 of these points are involutory points
(corresponding to the twelve lines of the corresponding
double-six, see the Theorem \ref{T:plane1}) and the other 15
belong respectively to other 15 Cremona  planes. It is natural to
expect that this configuration of $36$ planes can be obtained  by
cutting with a linear space the $36$ hyperplanes in $\P^5$
considered at 6.1.5.1 of \cite{Hu}, related to the Weyl group of
the exceptional group $E_6$, which have exactly the same
properties. \end{remark}

\begin{remark}\label{R:analytic}\rm The two components of the locus of singular L\"uroth quartics can
be also interpreted analytically, as follows. Consider the general
equation
$$\sum_{k=0}^4 \lambda_k\ell_0\cdots \hat\ell_k \cdots \ell_4 =0$$
of a L\"uroth quartic with inscribed pentalateral $\{\ell_0,
\dots, \ell_4\}$, considered in (21) of \cite{OS09}.
 Quartics in $\L_2$ have three of the five lines $\ell_k$ which are concurrent at the same point.
This follows easily from the identity:
$$
\left|\begin{array}{cccc}
\ell'_1+\ell'_0&\ell'_0&\ell'_0&\ell'_0\\
\ell'_0&\ell'_2+\ell'_0&\ell'_0&\ell'_0\\
\ell'_0&\ell'_0&\ell'_3+\ell'_0&\ell'_0\\
\ell'_0&\ell'_0&\ell'_0&\ell'_4+\ell'_0\\
\end{array}\right|=
\frac{1}{\lambda_0\lambda_1\lambda_2\lambda_3\lambda_4}
\sum_{k=0}^{4}\lambda_k\ell_0\cdots\hat{\ell_k}\cdots
\ell_{4}
$$
where $\ell'_k:=\frac{\ell_k}{\lambda_k}$.  Indeed, if $\ell_0,
\ell_1, \ell_2$ are concurrent at the same point, then the matrix
evaluated at this point has rank two. This condition corresponds
to the vanishing of the invariant $I$, as already remarked in \S
\ref{S:symm}.
 Quartics in $\L_1$ can be obtained for any
 given pentalateral, by a convenient choice of constants $\lambda_i$.
{\it Summarizing:  by specializing the $\ell_k$'s we  obtain
quartics in $\L_2$, and by specializing the $\lambda_k$'s we
obtain quartics in $\L_1$.}
\end{remark}

\emph{Proof of   Theorem \ref{T:main0}} We recall that
$s^*(\widetilde \D)$ consists of  the divisor of the twenty-seven
lines, with multiplicity two. The scheme $s^*(\widetilde \L_1)$ is
supported on the involutory points, and its length (on each line)
can be computed as the difference between the length of
$s^*(\widetilde \D\cap \widetilde \L)$ and the length of
$s^*(\widetilde \L_2)$ (both on a line); precisely it is equal to
$72-40=32$. Since the $16$ involutory points on each line are
distinct (see the proof of Prop. 6.3 in \cite{OS09}), it follows
that  $s^*(\widetilde \L_1)$ consists of $16$ points of length
$2$, in particular $s^*(\widetilde \L)$ is transversal to the
lines of $S$ at the involutory points, which implies   formula
(i). Then (ii) follows easily from the computation $27\cdot
54-27\cdot 24=27\cdot 30$.
  \qed

\begin{remark}\label{R:barth}\rm
From Theorem \ref{T:main0} it follows that the Barth map is not ramified over $\L_1$,
which answers  a question posed originally by Peskine, see
\cite{PT},  6.3 and 9.2.
\end{remark}

Part (ii) of the following theorem contains   Theorem \ref{noninvolutory}
of the introduction, with additional information.

\begin{theorem}\label{noninvolutory2}
Let $S$ be a nonsingular cubic surface.

\begin{itemize}

\item[(i)]
 Given a line  $\ell$   on  $S$ consider  the five planes $\Pi_i$, $i=1,\dots, 5$,
containing $\ell$ and such that $S\cap\Pi_i$  consists of $\ell$
and of two residual lines. Let $P_i$ be the intersection point of
the two residual lines. For every choice of distinct points $P_i,
P_j, P_k \in \{P_1, \dots, P_5\}$
  let $\Pi_{ijk}=<P_i, P_j, P_k>$ be the plane they span. Then
\[
\{\ell\cap \Pi_{ijk}: 1\le i < j < k \le 5\}
\]
are the ten non-involutory points on $\ell$.

\item[(ii)] Fix a double-six on $S$. Let $\ell_s$,  $s=1,\ldots
,15$, be the $15$  remaining lines. For each $1\le s \le 15$
consider the three planes $\Pi_{s,h}$,    $h=1,2,3$, containing
$\ell_s$ such that $S\cap\Pi_{s,h}$ consists of $\ell_s$ and of
two residual lines not belonging to the double-six. Let $P_{s,h}$
be the intersection point of the two residual lines. Then the
$15$ points $\ell_s\ \cap <P_{s,1}, P_{s,2}, P_{s,3}>$, $s=1,
\dots, 15$, lie on a plane, which is the Cremona plane associated
to the double-six. In particular the points $\ell_s\ \cap
<P_{s,1}, P_{s,2}, P_{s,3}>$, $s=1, \dots, 15$, are non-involutory
points on $\ell_s$.

\end{itemize}
\end{theorem}

\proof It is enough to prove (ii). Write the equation of $S$ in  Cremona form (\ref{E:cre1}) and
take   $\ell_s$ to be
$Z_0+Z_1=Z_2+Z_3=Z_4+Z_5=0$.
 The three pairs of lines coplanar with  $\ell_s$ are:
\[
\begin{array}{ll}
Z_0+Z_1=Z_2+Z_4=Z_3+Z_5=0,& Z_0+Z_1=Z_2+Z_5=Z_3+Z_4=0 \\
Z_2+Z_3=Z_0+Z_4=Z_1+Z_5=0,& Z_2+Z_3=Z_0+Z_5=Z_1+Z_4=0 \\
Z_4+Z_5=Z_0+Z_2=Z_1+Z_4=0, & Z_4+Z_5=Z_0+Z_4=Z_1+Z_2=0
\end{array}
\]
 The intersection points of the three pairs are respectively
{\small
\[
\begin{array}{ll}
P_{s,1}=\left(-( \beta_2+\beta_3-\beta_4-\beta_5),
\beta_2+\beta_3-\beta_4-\beta_5, \beta_0-\beta_1,\beta_0-\beta_1,-(\beta_0-\beta_1),-(\beta_0-\beta_1)\right) \\
P_{s,2}=\left(\beta_2-\beta_3,\beta_2-\beta_3,-(\beta_0+\beta_1-\beta_4-\beta_5),
\beta_0+\beta_1-\beta_4-\beta_5, -(\beta_2-\beta_3),-(\beta_2-\beta_3)\right) \\
P_{s,3}=\left(\beta_4-\beta_5,\beta_4-\beta_5,-(\beta_4-\beta_5),-(\beta_4-\beta_5),
-(\beta_0+\beta_1-\beta_2-\beta_3),\beta_0+\beta_1-\beta_2-\beta_3\right)
\end{array}
\]
}
The plane spanned by the three points cuts $\ell_s$ at the point
$$-(\beta_2-\beta_3)(\beta_4-\beta_5)P_{s,1}+
(\beta_0-\beta_1)(\beta_4-\beta_5)P_{s,2}-(\beta_0-\beta_1)(\beta_2-\beta_3)P_{s,3}=$$
{\tiny
$$=\left(\frac{\beta_2+\beta_3-\beta_4-\beta_5}{\beta_0-\beta_1},
-\frac{\beta_2+\beta_3-\beta_4-\beta_5}{\beta_0-\beta_1},
-\frac{\beta_0+\beta_1-\beta_4-\beta_5}{\beta_2-\beta_3},
\frac{\beta_0+\beta_1-\beta_4-\beta_5}{\beta_2-\beta_3},
\frac{\beta_0+\beta_1-\beta_2-\beta_3}{\beta_4-\beta_5},
-\frac{\beta_0+\beta_1-\beta_2-\beta_3}{\beta_4-\beta_5}\right)$$}
and a direct computation shows that this point belongs to the
Cremona plane
$\beta_0^2Z_0+\beta_1^2Z_1+\beta_2^2Z_2+\beta_3^2Z_3+\beta_4^2Z_4+\beta_5^2Z_5=0$.
Since this computation can be repeated for all the 15  lines not
belonging to the double six corresponding to the given Cremona
equations, we obtain the conclusion. \qed

There is a more conceptual proof of the above theorem, obtained by rephrasing an argument of Le Potier-Tikhomirov
in the geometry of the cubic surface,
and indeed we discovered it in this way.
 In \cite{PT}, Prop. 7.1,  they consider a nodal quartic with its six bitangent lines passing through the node.
They prove that the six contact points lie on a conic. If the
quartic is in $\L_2$ then this conic is singular (and also the
converse holds).
 On the   cubic surface $S$, when $P$ belongs to a line $\ell$ in $S$,
the ramification quartic of the projection  centered at $P$ has
the six bitangents passing through the node corresponding to the
five planes $\Pi_i$ and  to the tangent plane at $P$ itself. The
contact points of the six bitangents correspond to the lines
$\langle PP_i\rangle$ and the sixth one corresponds to the tangent
of the residual conic cut by the tangent plane at $P$. Hence these
six lines lie in a quadric cone with vertex at $P$. This quadric
cone splits into two planes when three contact points lie on a
line and this in turn corresponds to the fact that $P, P_i, P_j,
P_k$ lie on the same plane. So $P=\ell \ \cap <P_i, P_j, P_k>$ is
a non-involutory point, as in the proof given above.

We now  come to the question of  giving an explicit expression of
the L\"uroth invariant $\L$. We have used the description of
invariants of plane quartics given in \cite{Br}. From this
description we have shown, with a brute force computer
computation, that

\begin{prop}\label{noinv}
There is no invariant of degree $\le 15$ which vanishes on $(\L_2)_{red}$.
\end{prop}

The existence of such an invariant was asked by Morley in the last paragraph of \cite{fM14}.
This is interesting because from   Prop. \ref{noinv}  it follows that $(\L_2)_{red}$ is a divisor
 on the L\"uroth hypersurface
which is not cut by hypersurfaces.
 This makes it conceivable that
there is an invariant $I_{30}$ of degree $30$  which vanishes (doubly) on $(\L_2)_{red}$,
and that there is an invariant $I_{24}$ of degree $24$  which vanishes on $\L_1$.

 \noindent
Morley gives an explicit  form of the  restriction of the
invariant of degree 24 to a nodal cubic surface, assuming that it
exists. But he does not show that his formula  can be extended to
all cubic surfaces, and therefore he does not prove the existence
of the invariant. If this invariant exists, then $\L_1$ is a
complete intersection. This might be true, but it has not yet been
proved.

 These speculations are important
 because they give a hint about a possible explicit description of the  L\"uroth invariant.
Morley suggested the form
$\L=I_{27}\cdot I'_{27}+I_{24}I_{15}^2$
(up to scalar constants in the invariants appearing in the formula)
where $I_{27}=\D$ is the discriminant.
By   Prop. \ref{noinv}  this is not possible, but still the expression
$\L=I_{27}\cdot I'_{27}+I_{24}I_{30}$
could be possible, and we leave  its existence  as a question.

\noindent
 This expression means that, modulo the discriminant,
the L\"uroth invariant is the product of the two invariants $I_{24}$ and $I_{30}$.
The computation of $\L$  could reduce in this case to smaller degree invariants.
Let us recall that the main result of Dixmier \cite{Dix} states that
every invariant of a plane quartics is algebraically
dependent on an explicit system of invariants of degrees up to $27$.

\noindent Moreover, Morley asks if $I'_{27}$ is the discriminant
too. This can be probably answered by a computational analysis,
which does not seem easy and we do not pursue here.

The following result says something about the singularities of the
L\"uroth hypersurface $\mathcal{L}$:

\begin{prop}\label{P:nnormal}
The L\"uroth hypersurface $\mathcal{L}$ is not normal.
 \end{prop}

 \proof
 Consider the incidence relation:
   \[
   \widetilde{\L} := \left\{(B,\{\ell_0,\dots,\ell_4\}):
   \begin{array}{l}
   \hbox{$\{\ell_0,\dots,\ell_4\}$ is a complete} \\
\hbox{ 5-lateral and $B$ is a n.s. }\\
   \hbox{ quartic circumscribed to it}\end{array}\right\}  \subset  \P^{14}\times (\P^{2\vee})^{(5)}
   \]
   and the projections:
   \[
   \xymatrix{
   &  \widetilde{\L} \ar[dl]_{q_1}\ar[dr]^{q_2} \\
   \P^{14}&&(\P^{2\vee})^{(5)} }
   \]
   Clearly   $\overline{q_1(\widetilde{\L})}=\mathcal{L} \subset \P^{14}$.
     $\widetilde{\L}$ is irreducible of
   dimension 14, and  its nonsingular locus contains
   $q_2^{-1}(U)$, where $U \subset (\P^{2\vee})^{(5)}$ is the
   locus of strict pentalaterals (i.e. those having 10 distinct
   vertices). The  fibre $q_1^{-1}(C)$ over a general $[C]\in \mathcal{L}$
    is one-dimensional and irreducible, consisting of the 5-laterals which are
    inscribed in $C$. But  there is a class of  nonsingular L\"uroth
   quartics, the  \emph{desmic quartics}, such that $q_1^{-1}(C)$
   is disconnected and each one of its component
   intersects $q_2^{-1}(U)$ (see for example \cite{hB14}).
   Consider the morphism
   \[
\xymatrix{q_2^{-1}(U)\ar[r] & \P^{14} }
\]
and its Stein factorization:
\[
\xymatrix{q_2^{-1}(U)\ar[r] &\mathcal{S}\ar[r]^-\upsilon &\P^{14}
}
\]
From the above description it follows that $\upsilon$ maps
$\mathcal{S}$ birationally and dominantly to $\mathcal{L}$, but it
has disconnected fibres over the locus of desmic quartics. From
Zariski's Main Theorem it follows that $\mathcal{L}$ is not normal
along the locus of desmic quartics. In particular the singular
locus of $\mathcal{L}$ has codimension $1$ in $\mathcal{L}$. \qed

\section{The class of the   L\"uroth divisor in $\overline{M}_3$}

Let $\overline{M}_3$ be the coarse  moduli space  of stable curves
of genus 3. We will compute certain rational divisor classes in
$\mathrm{Pic}_{\mathrm{fun}}(\overline{M}_3)\otimes \mathbb{Q}$ in
terms of the Hodge class $\lambda$ and  of $\delta_0,\delta_1$,
the classes of the boundary components $\Delta_0,\Delta_1$.
Precisely, $\Delta_0$  generically parametrizes irreducible
 singular stable curves, and $\Delta_1$ generically
parametrizes reducible stable curves.

The hyperelliptic locus $\overline{H} \subset \overline{M}_3$ is a
divisor whose class is
\[
\overline{h} = 9\lambda - \delta_0-3\delta_1
\]
  For the
proof see \cite{HM98} and \cite{mT88}.

The discriminant $\mathcal{D}\subset \P^{14}$, being an
SL$(3)$-invariant hypersurface of degree 27, admits a rational map
  $\Phi:\xymatrix{\mathcal{D}\ar@{-->}[r]&\overline{M}_3}$.
  The
closure of the image is a divisor
$D:=\overline{\mathrm{Im}(\Phi)}\subset \overline{M}_3$. Since
$\mathcal{D}$ contains the double conics and the cuspidal curves,
$D$ contains $\overline{H}$, $\Delta_0$ and $\Delta_1$.
 In \cite{AC93} it is proved that
$\mathcal{D}$ vanishes with multiplicity 14 on double conics and
with multiplicity 2 on cuspidal quartics. If we consider any test
curve intersecting double conics and cuspidal quartics
transversely, we will have to perform a base change in order to
get stable reduction. Standard facts about stable reduction (see
\cite{HM98}, Chapter 3.C)    imply that the degrees of the base
changes needed are respectively 2 and 6. It follows that  $[D]$
contains $\overline{H}$ and $\Delta_1$ with multiplicity $2\cdot
14$ and $6\cdot 2$ respectively.   Therefore the class $[D]$ is
computed as follows:
\[
\begin{array}{ll}
[D] &= 2\cdot 14 [H]+6\cdot 2 \delta_1 + \delta_0 \\
&= 28 (9\lambda - \delta_0-3\delta_1) + 12 \delta_1
+ \delta_0 \\
&= 9(28\lambda-3\delta_0-8\delta_1)
\end{array}
\]
This formula can be easily tested and confirmed using the pencils
considered in Exercise (3.166) of \cite{HM98}.

 We define \emph{the
divisor $L\subset \overline{M}_3$ of Luroth quartics}, to be  the
closure of the locus of nonsingular Luroth quartics.

\begin{prop}\label{P:lurvanish}
Let $\mathcal{L}\subset \P^{14}$ be the hypersurface
  of Luroth quartics. Then $\mathcal{L}$ does not contain the
  loci of double conics and of cuspidal quartics.   Moreover we
  have:
  \[
[L] = 18(28\lambda-3\delta_0-8\delta_1)
\]
  \end{prop}

\proof       From Proposition \ref{geominvolutory} it follows that
general quartics in $\mathcal{L}_1$ and in $\mathcal{L}_2$ can be
obtained as branch curves of projections from points belonging to
a line $\ell$   of a general cubic surface $S$. It is well known
that such branch curves have one node, which is the projection of
the line $\ell$, with principal tangents the projections of the
two planes containing  $\ell$  and such that the residual conic is
tangent to $\ell$. Therefore  both $\mathcal{L}_1$ and
$\mathcal{L}_2$ contain a nodal quartic.  Since $\mathcal{L} \cap
\mathcal{D}= \mathcal{L}_1\cup \mathcal{L}_2$ has pure dimension
12, as well as the locus of cuspidal quartics, it follows that
this locus is not entirely contained in $\mathcal{L}$.

For double conics we argue as follows. Consider an irreducible
quartic $C$ with a node $O$ but otherwise nonsingular. Then the
six contact points of $C$ with the tangents from $O$ (which are
the intersections different from $O$ of $C$ with the first polar
of $O$ with respect to $C$) are on a
  conic $\theta$, and this conic is reducible
   if and only if $C \in  \mathcal{L}_2$
  (\cite{PT}, \S 7).
   If the quartic $C$   degenerates to a double conic $2\vartheta$
   then $\vartheta$ is a component of the first polar of any of
   its points. This implies that $\vartheta$ is a degeneration of
   $\theta$.
   Then,
  since $\theta$ is reducible if $C\in  \mathcal{L}_2$,
  it follows that $C$ cannot degenerate to a double
   irreducible conic. Therefore $\mathcal{L}_2$ does not contain the locus of
  double conics.
 Consider now a
quartic  $C\in\mathcal{L}_1$, such that
$C\notin\mathcal{L}_1\cap\mathcal{L}_2$. Then $C$ is circumscribed
to a 5-lateral having 10 distinct vertices (Remark
\ref{R:analytic}). Therefore $C$ cannot be a double irreducible
conic.  In conclusion,  $\mathcal{L}$ does not contain double
irreducible conics.

  Being $\mathrm{SL}(3)$-invariant, $\mathcal{L}$ defines a divisor in
$\overline{M}_3$, containing $L$, whose class is
\[
[\mathcal{L}] = 2[D]
\]
because $\mathrm{deg}(\mathcal{L})=54= 2\
\mathrm{deg}(\mathcal{D})$. On the other hand, we have
\[
[L] = [\mathcal{L}]- a  \overline{h} - b\delta_1 = 2 [D]- a
\overline{h} - b\delta_1
\]
where $a,b$ are the multiplicities of $\mathcal{L}$ along the loci
of double conics and of cuspidal quartics respectively.  From the
first part of the proof it follows that $a=b=0$. Therefore
\[
[L] = 2[D] = 18(28\lambda-3\delta_0-8\delta_1)
\]
 \qed

 In a similar vein one  computes the class $[\mathbf{Cat}]$ of
the catalecticant hypersurface to be:
\[
[\mathbf{Cat}]  = 2(28\lambda-3\delta_0-8\delta_1)
\]

Note that $28\lambda-3\delta_0-8\delta_1$ is necessarily the class
of the divisor defined by the SL$(3)$ invariant of (smallest)
degree 3.

  \vskip1cm

 \noindent
  \textsc{g. ottaviani} -
  Dipartimento di Matematica ``U. Dini'',
  Universit\`a di Firenze, viale Morgagni 67/A, 50134 Firenze (Italy). e-mail: \texttt{ottavian@math.unifi.it}

  \separation
  \noindent
\textsc{e. sernesi} -
 Dipartimento di Matematica,
 Universit\`a Roma Tre,
 Largo S.L. Murialdo 1,
 00146 Roma (Italy).  e-mail: \texttt{sernesi@mat.uniroma3.it}

\end{document}